\pdfoutput=1
\RequirePackage{ifpdf}
\ifpdf 
\documentclass[pdftex]{sigma}
\else
\documentclass{sigma}
\fi

\usepackage{tikz}
\usetikzlibrary{intersections,calc,positioning}
\usetikzlibrary{shapes.geometric,fit}
\usetikzlibrary{decorations,arrows,decorations.pathmorphing,backgrounds,positioning,fit,petri,shadows}

\newcommand{\Z}
{{\mathbb Z}}

\newcommand{\g}
{{\mathfrak g}}
\newcommand{\lgg}
{{\mathfrak l}}

\newcommand{\lggt}
{\tilde{\mathfrak l}}
\newcommand{\gk}
{\hat{{\mathfrak g}}}
\newcommand{\gt}
{\tilde{{\mathfrak g}}}

\newcommand{\h}
{{\mathfrak h}}

\newcommand{\ox}
{{\overline{x}}}

\numberwithin{equation}{section}

\begin{document}
\allowdisplaybreaks

\newcommand{\arXivNumber}{2504.15597}

\renewcommand{\thefootnote}{}

\renewcommand{\PaperNumber}{071}

\FirstPageHeading

\ShortArticleName{Linear Independence for $A_1^{(1)}$ by Using $C_{2}^{(1)}$}

\ArticleName{Linear Independence for $\boldsymbol{A_1^{(1)}}$ by Using $\boldsymbol{C_{2}^{(1)}}$\footnote{This paper is a~contribution to the Special Issue on Recent Advances in Vertex Operator Algebras in honor of James Lepowsky. The~full collection is available at \href{https://www.emis.de/journals/SIGMA/Lepowsky.html}{https://www.emis.de/journals/SIGMA/Lepowsky.html}}}

\Author{Mirko PRIMC~$^{\rm a}$ and Goran TRUP\v{C}EVI\'C~$^{\rm b}$}

\AuthorNameForHeading{M.~Primc and G. Trup\v {c}evi\'c}

\Address{$^{\rm a)}$~Faculty of Science, University of Zagreb, Zagreb, Croatia}
\EmailD{\href{mailto:primc@math.hr}{primc@math.hr}}

\Address{$^{\rm b)}$~Faculty of Teacher Education, University of Zagreb, Zagreb, Croatia}
\EmailD{\href{mailto:goran.trupcevic@ufzg.hr}{goran.trupcevic@ufzg.hr}}
	
\ArticleDates{Received May 02, 2025, in final form August 12, 2025; Published online August 19, 2025}
	
\Abstract{In the previous paper, the authors proved linear independence of the combinatorial spanning set for standard \smash{$C_\ell^{(1)}$}-module $L(k\Lambda_0)$ by establishing a connection with the combinatorial basis of Feigin--Stoyanovsky's type subspace $W(k\Lambda_0)$ of \smash{$C_{2\ell}^{(1)}$}-module $L(k\Lambda_0)$. In this note we extend this argument for \smash{$C_{1}^{(1)}\cong A_{1}^{(1)}$} to all standard \smash{$A_{1}^{(1)}$}-modules $L(\Lambda)$. In the proof we use a coefficient of an intertwining operator of the type \smash{$\binom{L(\Lambda_2)}{L(\Lambda_1)\ L(\Lambda_1)}$} for standard \smash{$C_{2}^{(1)}$}-modules.	}
	
\Keywords{affine Lie algebras; standard modules; Feigin--Stoyanovsky's type subspace; combinatorial basis}
	
\Classification{17B67; 17B69}

\renewcommand{\thefootnote}{\arabic{footnote}}
\setcounter{footnote}{0}

\section{Introduction}

In \cite{PT}, the authors proved linear independence of the combinatorial spanning set for standard \smash{$C_\ell^{(1)}$}-module $L(k\Lambda_0)$ by establishing a connection with the combinatorial basis of Feigin--Stoyanovsky's type subspace $W(k\Lambda_0)$ of \smash{$C_{2\ell}^{(1)}$}-module $L(k\Lambda_0)$ constructed in \cite{BPT}. For $\ell=1$ we have \smash{$C_{1}^{(1)}\cong A_{1}^{(1)}$} and the combinatorial basis of \smash{$L_{A_{1}^{(1)}}(k\Lambda_0)$} is a part of the general construction of combinatorial bases of all standard \smash{$A_{1}^{(1)}$}-modules \smash{$L_{A_{1}^{(1)}}(k_0\Lambda_0+k_1\Lambda_1)$}, $k_0+k_1=k$, obtained independently in \cite{MP} and \cite{FKLMM}.
In this note we extend the argument from \cite{PT} to all standard \smash{$A_{1}^{(1)}$}-modules by using a coefficient of an intertwining operator of the type \smash{$\binom{L(\Lambda_2)}{L(\Lambda_1)\ L(\Lambda_1)}$} for standard \smash{$C_{2}^{(1)}$}-modules. This gives a new proof of linear independence of combinatorial bases of standard \smash{$A_{1}^{(1)}$}-modules and, hopefully, this approach may lead to a proof of linear independence of combinatorial bases of all standard \smash{$C_{\ell}^{(1)}$}-modules conjectured in \cite{CMPP}.

As in \cite{PT}, the key idea for the proof of linear independence of the spanning set $B_1$ of monomial vectors $x(\pi)v_{\bar\Lambda}$ in \smash{$L_{A_{1}^{(1)}}(\bar\Lambda)$} is to embed Lie algebra $\lgg$ of type $A_1$ into ${\mathfrak g}$ of type $C_{2}$ (see Figure \ref{C2}) together with its standard module
\[
L_{A_{1}^{(1)}}\bigl(\bar{\Lambda}\bigr)\subset L_{C_{2}^{(1)}}(\Lambda)
\supset W_{C_{2}^{(1)}}(\Lambda),
\]
and then, by using the inner derivations $T$ of ${\mathfrak g}$, connect the set $B_1$ and
the basis $B_2^1$ of the Feigin--Stoyanovsky subspace \smash{$W_{C_{2}^{(1)}}(\Lambda)$} consisting of monomial vectors $\ox(\pi)v_{\Lambda}$.

Monomials $x(\pi)$ and $\ox(\pi)$ in the universal enveloping algebra $U(\hat{\mathfrak g})$, given by \eqref{upi} and~\eqref{wpi}, are parameterized with colored partitions $\pi$ in three colors, determined by frequencies $\{a_j, b_j, c_j\mid j\geq0\}$ satisfying {\it the same} difference conditions \eqref{DCA1}--\eqref{DCA2}. The case when $\bar\Lambda=k\Lambda_0$ and $\Lambda=k\Lambda_0$ is relatively simple because $\pi$ in monomials $x(\pi)$ and $\ox(\pi)$ satisfy the {\it the same} initial conditions~${a_0=b_0=c_0=0}$ and for the proper power $T^{N'}$ of $T$ we have
\[
T^{N'}\colon\ x(\pi)\to \ox(\pi).
\]
However, in general the initial conditions for $\pi$ in monomials $x(\pi)$ and $\ox(\pi)$ {\it are not} the same, see \eqref{ICA1}--\eqref{ICA2} compared to \eqref{ICC1}--\eqref{ICC2}, and, together with $T$, a~coefficient $w$ of an intertwining operator of the type \smash{$\binom{L(\Lambda_2)}{L(\Lambda_1)\ L(\Lambda_1)}$} for standard \smash{$C_{2}^{(1)}$}-modules is used to circumvent this difficulty.

\section[Affine Lie algebras of type A\_1\^(1) C\_2\^(1) and their standard modules]{Affine Lie algebras of type $\boldsymbol{ A_1^{(1)} \subset C_2^{(1)}}$\\ and their standard modules}

\subsection{Affine Lie algebras and standard modules}

Let $\g$ be a simple Lie algebra with a Cartan decomposition $\mathfrak g={\mathfrak h}+\sum_{\alpha\in R} \mathfrak g_\alpha$, where
$R$
 is a root system of $\g$. Let $\alpha_1,\dots, \alpha_n$ be a basis of $R$, $\theta$ the maximal root, and $\omega_1, \dots,\omega_n$ the corresponding fundamental weights. For each root $\alpha$ fix a root vector $x_\alpha$ in $\g_\alpha$.
Let $\langle\ ,\ \rangle$ be the normalized Killing form so that $\langle\theta,\theta\rangle=2$, through which we identify $\h$ and $\h^*$.
Let $B$ be a basis of ${\mathfrak g}$ consisting of root vectors and elements of $\h$, with a linear order $\succ$.

Let $\tilde{\mathfrak g}$ be the affine Lie algebra associated to $\g$,
\[
\hat{\mathfrak g} = \mathfrak{g} \otimes
\mathbb{C}\bigl[t,t^{-1}\bigr] + \mathbb{C}c, \qquad
\tilde{\mathfrak g} = \hat{\mathfrak g} + \mathbb{C} d,
\]
with commutation relations
\[
[x(i),y(j)]= [x,y](i+j)+ i\langle x,y\rangle
\delta_{i+j,0}c,\qquad [c,\gt] = 0, \qquad [d,x(j)] = j x(j),\]
where $x(n)=x\otimes t^{n}$, for $x\in{\mathfrak g}$ and $n\in\mathbb Z$. Identify $\g=\g\otimes 1\subset \gk$.
Set $\overline{B}=\{b(n) \mid b\in B,\, n \in\mathbb{Z}\}$; so that~${\overline{B}\cup \{c\}}$ is a basis of $\gk$. Extend the order on $B$ to $\overline{B}$:
$b(n)\succ b'(n')$ if $n>n'$ or $n=n'$, $ b\succ b'$.

Denote by $\Lambda_0,\dots,\Lambda_n$ the fundamental weights of $\gt$. For a given
$\Lambda=k_0 \Lambda_0+k_1 \Lambda_1+\dots+k_n \Lambda_n$, let $L(\Lambda)=L_{\gt}(\Lambda)$ be a standard (i.e., integrable highest weight) $\tilde{\mathfrak g}$-module, $v_\Lambda$ a fixed highest weight vector, and $k=\Lambda(c)$ the level of $L(\Lambda)$ (cf.\ \cite{K}).


\subsection[Bases of standard modules for affine Lie algebra of type A\_1\^(1)subset C\_2\^(1)]{Bases of standard modules for affine Lie algebra of type $\boldsymbol{ A_1^{(1)}\subset C_2^{(1)}}$}

Let $\g$ be a simple Lie algebra of the type $C_2$. Let
\begin{equation} \label{E: root system}
	R=\{2\epsilon_1, \epsilon_1 + \epsilon_2, 2\epsilon_2, \epsilon_1 - \epsilon_2, \epsilon_2 - \epsilon_1, - 2 \epsilon_2, -\epsilon_1 - \epsilon_2,-2\epsilon_1\} \subset \mathbb R^2
\end{equation}
be a root system of $\g$. Let $\alpha_1=\epsilon_1 - \epsilon_2$, $\alpha_2=2 \epsilon_2$ be a root basis, $\theta=2\epsilon_1$ the maximal root, and~${\omega_1=\epsilon_1}$, $\omega_2=\epsilon_1+\epsilon_2$ the corresponding fundamental weights.

Fix root vectors
$x_{11}$, $x_{12}$, $x_{22}$, $x_{1 \underline{2}}$, $x_{2\underline{1}}$, $x_{\underline{2}\underline{2}}$, $x_{\underline{2}\underline{1}}$, $x_{\underline{1}\underline{1}}$
corresponding respectively to the roots in~\eqref{E: root system}
and let $x_{1\underline{1}}$, $x_{2 \underline{2}}$
be the simple coroots in $\mathfrak h$ corresponding to positive roots $2\epsilon_1$ and~$2\epsilon_2$.
These vectors form a weight basis $B$ of ${\mathfrak g}$. Define an order on $B$ in the following way: set~${1\succ 2 \succ \underline{2} \succ \underline {1}}$
and define a lexicographic order
$x_{ab}\succ x_{a'b'}$ if $a\succ a'$ or $a=a'$, $b\succ b'$.

Denote by $\Lambda_0$, $\Lambda_1$, $\Lambda_2$ the fundamental weights of $\gt$. For
$\Lambda=k_0 \Lambda_0+k_1 \Lambda_1+k_2 \Lambda_2$, let \smash{$L(\Lambda)=L_{C_2^{(1)}}(\Lambda)$} be a standard $\tilde{\mathfrak g}$-module of the level $k=k_0+k_1+k_2$ with a fixed highest weight vector $v_\Lambda$. The standard module $L(\Lambda)$ can be realised in the tensor product of level $1$ standard modules,
$L(\Lambda)\subset L(\Lambda_0)^{\otimes k_0}\otimes L(\Lambda_1)^{\otimes k_1} \otimes L(\Lambda_2)^{\otimes k_2}$.

On the top of the standard module $L(\Lambda_0)$ there is the trivial $1$-dimensional module for $\g$; denote by $v_0$ its weight vector. On the top of $L(\Lambda_1)$ there is the $4$-dimensional irreducible module~$V(\omega_1)$ for $\g$ with weights $\epsilon_1$, $\epsilon_2$, $-\epsilon_1$, $-\epsilon_2$ and corresponding weight vectors $v_1$, $v_2$, $v_{\underline{1}}$, $v_{\underline{2}}$. On the top of $L(\Lambda_2)$ there is the $5$-dimensional irreducible module $V(\omega_2)$ for $\g$ with weights $\epsilon_1+\epsilon_2$, $-\epsilon_1+\epsilon_2$, $0$, $\epsilon_1-\epsilon_2$, $-\epsilon_1-\epsilon_2$ and corresponding weight vectors $v_{12}$, $v_{\underline{1}2}$, $v_{00}$, $v_{1 \underline{2}}$, $v_{\underline{1}\underline{2}}$. Note that $v_0=v_{\Lambda_0}$, $v_1=v_{\Lambda_1}$, $v_{12}=v_{\Lambda_2}$.

The subalgebra $\lgg=\operatorname{span}\{x_{11},x_{1 \underline{1}}, x_{\underline{1}\underline{1}}\}\subset\g$
is a simple algebra of type $A_1$ with the simple root $\theta=2\epsilon_1$.

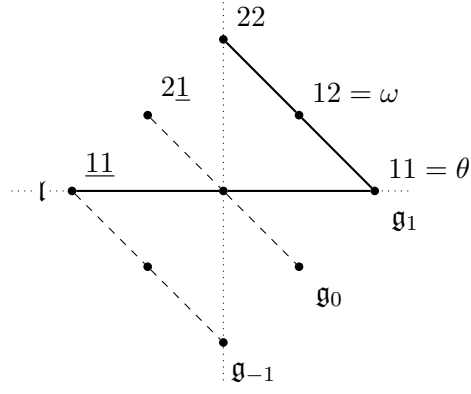
\begin{figure}[ht]\centering
		\begin{tikzpicture}
			
			\draw[dotted] (-2.8,0) -- (2.5,0);
			\draw[dotted] (0,-2.5) -- (0,2.5);,
			
			\node[draw,fill, circle, inner sep=1pt,label=45:${11=\theta}$](11) at (2,0) {};
			\node[draw,fill, circle, inner sep=1pt,label=45:${12=\omega}$](12) at (1,1) {};
			\node[draw,fill, circle, inner sep=1pt,label=45:$22$](22) at (0,2) {};
			\node[draw,fill, circle, inner sep=1pt](1b) at (1,-1) {};
			\node[draw,fill, circle, inner sep=1pt](0) at (0,0) {};
			\node[draw,fill, circle, inner sep=1pt,label=45:$2\underline{1}$](2a) at (-1,1) {};
			\node[draw,fill, circle, inner sep=1pt](ba) at (-1,-1) {};
			\node[draw,fill, circle, inner sep=1pt,label=45:$\underline{1}\underline{1}$](aa) at (-2,0) {};
			\node[draw,fill, circle, inner sep=1pt](bb) at (0,-2) {};
			\draw[
			thick] (11) -- (12) -- (22) ;
			\draw[
			thick] (11) -- (0) -- (aa);
			\draw[dashed] (aa) -- (ba) -- (bb);
			\draw[dashed] (2a) -- (0) -- (1b);
			
			\node[fill=white, inner sep=1pt] at (-2.4,0)	{$\lgg$};				
			\node at (2.4,-.4)	{$\g_1$};				
			\node at (1.4,-1.4)	{$\g_0$};				
			\node at (.4,-2.4)	{$\g_{-1}$};				
		\end{tikzpicture}
 \caption{The root system of type $C_2$.} \label{C2}
\end{figure}

The inclusion $\lgg\subset\g$ induces an inclusion of affine Lie algebras $\lggt\subset\gt$; the subalgebra $\lggt$ is of type \smash{$A_1^{(1)}$}. Denote by $\bar{\Lambda}_0$ and $\bar{\Lambda}_1$ fundamental weights of $\lggt$. Standard $\lggt$-modules can be found as $\lggt$-submodules of standard $\gt$-modules
\begin{gather*}
	L_{A_1^{(1)}} \bigl(\bar{\Lambda}_0\bigr)  \cong U\bigl(\lggt\bigr)v_{\Lambda_0} \subset L_{C_2^{(1)}} (\Lambda_0), \\	
	L_{A_1^{(1)}} \bigl(\bar{\Lambda}_1\bigr)  \cong U\bigl(\lggt\bigr)v_{\Lambda_1} \subset L_{C_2^{(1)}} (\Lambda_1)
	 \cong U\bigl(\lggt\bigr)v_{\Lambda_2} \subset L_{C_2^{(1)}} (\Lambda_2), \\ 	
	L_{A_1^{(1)}} \bigl(\bar{\Lambda}\bigr)  \cong U\bigl(\lggt\bigr)v_{\Lambda} \subset L_{C_2^{(1)}}(\Lambda), 	
\end{gather*}
for $\bar{\Lambda} = \bar{k}_0 \bar{\Lambda}_0 + \bar{k}_1 \bar{\Lambda}_1$,
$\Lambda = k_0 \Lambda_0 + k_1 \Lambda_1+ k_2 \Lambda_2$, where $\bar{k}_0=k_0$,
 $\bar{k}_1=k_1+k_2$.
For this reason we will use the notation
$\Lambda_0$, $\Lambda_1$, $\Lambda_2$ also for $\lggt$.

The PBW spanning set of \smash{$L_{A_1^{(1)}}(\Lambda)$},
$\Lambda = k_0 \Lambda_0 + k_1 \Lambda_1$,
 can be reduced to a monomial basis consisting of monomial vectors
\begin{equation} \label{upi}
	\prod_{j\geq 0} x_{\underline{1}\underline{1}}(-j)^{c_j} x_{1 \underline{1}}(-j)^{b_j} x_{11} (-j)^{a_j} v_\Lambda
\end{equation}
satisfying {\em difference conditions}
\begin{gather}
	\label{DCA1} 	
	a_i +b_i+a_{i+1}\leq k ,	\qquad
	c_i +b_i+a_{i+1}\leq k 	,\\
	\label{DCA1.5} c_i +b_{i+1}+a_{i+1}\leq k ,	\\
	\label{DCA2}	c_i +b_{i+1}+c_{i+1}\leq k 	
\end{gather}
and {\em initial conditions}
\begin{gather}
a_0=b_0  = 0,\qquad
	c_0 \leq k_1	,\label{ICA1} \\
a_1 \leq k_0	\label{ICA2}
\end{gather}
for \smash{$L_{A_1^{(1)}}(\Lambda)$}.
 This was proved by different methods in \cite{MP} and \cite{FKLMM}. In this note we give a~new proof of linear independence by transforming its elements to elements of a monomial basis of a~Feigin--Stoyanovsky subspace for $\gt$.


\subsection{Feigin--Stoyanovsky subspace}

For a simple Lie algebra $\g$, let $\omega$ be a minuscule coweight of $\g$, $\langle \omega, R\rangle = \{-1,0,1\}$. The minuscule coweight $\omega$ induces a $\Z$-gradation on $\g$
\[
\mathfrak g =\mathfrak g_{-1} + \mathfrak g_0 + \mathfrak g_1 , \qquad
{\mathfrak g}_0 = {\mathfrak h} +
\sum_{\langle\omega,\alpha\rangle=0}  {\mathfrak g}_\alpha,
\qquad
\displaystyle {\mathfrak g}_{\pm1} =
\sum_{\langle\omega,\alpha\rangle= \pm 1}  {\mathfrak g}_\alpha.
\]
The set
$
	\Gamma =
	\{ \alpha \in R \mid \langle\omega,\alpha\rangle = 1\}
$
is called the set of colors.
The subalgebras ${\mathfrak g}_{\pm1}\subset\mathfrak g$ are commutative, while $\g_0$ is reductive.

The $\Z$-gradation of $\g$ induces a $\Z$-gradation of $\gt$,
$
\gt =\gt_{-1} + \gt_0 + \gt_1$,
where
\[
\gt_0 = \g_0 \otimes \mathbb{C}\bigl[t,t^{-1}\bigr] + \mathbb{C}c + \mathbb{C} d,\qquad
\gt_{\pm 1} = \g_{\pm 1} \otimes \mathbb{C}\bigl[t,t^{-1}\bigr].\]
Again, the subalgebras $\gt_{\pm1}\subset\gt$ are commutative.

{\em Feigin--Stoyanovsky subspace} $W(\Lambda)=W_{\gt}(\Lambda)$ of a standard module $L_{\gt}(\Lambda)$ is a $\gt_1$-submodule generated by the highest weight vector
$W(\Lambda)=U(\gt_1)v_\Lambda$.


\subsection[Bases of Feigin--Stoyanovsky subspaces for affine Lie algebra of type C\_2\^(1)]{Bases of Feigin--Stoyanovsky subspaces for affine Lie algebra of type $\boldsymbol{ C_2^{(1)}}$}

Let $\g$ be a simple Lie algebra of type $C_2$, as before.
The minuscule coweight $\omega=\omega_2$ induces a~$\Z$-gradation on $\g$ with the set of colors $\Gamma = \{2\epsilon_1, \epsilon_1 + \epsilon_2, 2\epsilon_2\}$ (see Figure \ref{C2}).

The Feigin--Stoyanovsky subspace \smash{$W_{C_2^{(1)}}(\Lambda)$} for $ \Lambda=k_0\Lambda_0+k_1\Lambda_1+k_2\Lambda_2$, $k=k_0+k_1+k_2$, has a monomial basis
\begin{equation} \label{wpi}
\prod_{j\geq 0} x_{2 2}(-j)^{c_j} x_{1 2}(-j)^{b_j} x_{1 1}(-j)^{a_j} v_\Lambda
\end{equation}
satisfying
{\em difference conditions} \eqref{DCA1}--\eqref{DCA2}
and {\em initial conditions}
\begin{gather}
	\label{ICC1} a_0=b_0 =c_0  = 0,\qquad
	a_1 \leq k_0,\qquad
	a_1+b_1\leq k-k_2, \\
	\label{ICC2} b_1+c_1\leq k-k_2
\end{gather}
for \smash{$W_{C_2^{(1)}}(\Lambda)$} (cf.\ \cite{BPT} and \cite{P}\footnote{When interchanging $C_2^{(1)}\leftrightarrow B_2^{(1)}$ we should interchange $\Lambda_1 \leftrightarrow \Lambda_2$.}).


\section[Proof of linear independence for A\_1\^(1)]{Proof of linear independence for $\boldsymbol{ A_1^{(1)}}$}


\subsection[Translation L\_A\_1\^(1)(Lambda)to W\_C\_2\^(1)(Lambda')]{Translation $\boldsymbol{ L_{A_1^{(1)}}(\Lambda)\to W_{C_2^{(1)}}(\Lambda')}$}

Let $\Lambda = k_0 \Lambda_0 + k_1 \Lambda_1$. Let
\[\pi=\prod_{j\geq 0}(-j)^{c_j}(-j)^{b_j}(-j)^{a_j}\]
be a colored partition with parts $-j$ of colors $a$, $b$, $c$ appearing with frequencies $a_j$, $b_j$, $c_j$, which satisfy difference conditions \eqref{DCA1}--\eqref{DCA2}.
Denote by $x(\pi)$ and $\ox(\pi)$ monomials \eqref{upi} and \eqref{wpi}, respectively. Note that $x(\pi)$ is noncommutative and $\ox(\pi)$ is commutative.

Denote by $T=\textrm{ad}  x_{12}$ a derivative on $\gk$ and $U(\gk)$. Then, up to a scalar,
$ T x_{1 \underline{1}}(j)  = x_{1 2}(j)$,
 $T x_{\underline{1} \underline{1}}(j)  = x_{2\underline{1}}(j)$, $
 T x_{2\underline{1}}(j)  = x_{22}(j)$.
Note that
$T x_{1 1}(j) = T x_{1 2}(j)=0$.

For a monomial $x(\pi)$ satisfying difference and initial conditions for \smash{$L_{A_1^{(1)}}(\Lambda)$} set
\[N'=\sum_{j\geq 0} b_j + 2\sum_{j\geq 0} c_j.\]
The action by $T^{N'}$ transforms $x(\pi)$ to $\ox(\pi)$:
$T^{N'} x(\pi) = \ox (\pi)$.
Furthermore,
$T^{N'+1} x(\pi) = 0$.

Let $v_\Lambda=v_0^{\otimes k_0} \otimes v_1^{\otimes k_1}$ be a highest weight vector of \smash{$L_{C_2^{(1)}}(\Lambda)$}.
Then
\begin{align*}
	x_{12}(0)^{N'} x(\pi) v_\Lambda & = x_{12}(0)^{N'-1} (T x(\pi)) v_\Lambda + x_{12}(0)^{N'-1} x(\pi) x_{12}(0) v_\Lambda \\
	& = \dots= \bigl(T^{N'} x(\pi)\bigr) v_\Lambda= \ox(\pi) v_\Lambda.
\end{align*}
Note that $x_{12}(0)$ annihilates $v_\Lambda$ since the corresponding root is positive.
Furthermore, initial conditions for \smash{$W_{C_2^{(1)}}(\Lambda)$} imply that
$\ox(\pi) v_\Lambda \neq 0 \Leftrightarrow c_0=0$.

Hence for a monomial $x(\pi)$ set
\[N = N(\pi) =\sum_{j\geq 0} b_j + 2\sum_{j\geq 0} c_j - c_0.\]
Let $x(\pi)=x(\pi_1)x_{\underline{1} \underline{1}}(0)^{c_0}$. Note that $N(\pi) = N(\pi_1) + c_0$.
Then
\begin{align}
x_{12}(0)^N x(\pi) v_\Lambda & =  \bigl(T^N x(\pi)\bigr) v_\Lambda = \bigl(T^N \bigl(x(\pi_1)x_{\underline{1}\underline{1}}(0)^{c_0}\bigr)\bigr) v_\Lambda\nonumber\\
\label{eq2}	& = \ox(\pi_1)x_{2\underline{1}}(0)^{c_0} v_\Lambda + \sum \cdots x_{2 2}(0) v_\Lambda \\
\label{eq3}	& = \ox(\pi_1)x_{2\underline{1}}(0)^{c_0} v_\Lambda.
\end{align}

In \eqref{eq2} the sum goes over all the other possibilities of action of $T^N$ on factors in $x(\pi)$. In all of these at least two $T$'s act on the same $x_{\underline{1} \underline{1}}(0)$ factor. Hence one gets at least one $x_{2 2}(0)$ factor, which commutes with $x_{\underline{1} \underline{1}}(0)$ and $x_{2 \underline{1}}(0)$, and annihilates $v_\Lambda$.

Notice that $x_{2\underline{1}}(0)$ acts on the $v_1$'s in the tensor product \smash{$v_0^{\otimes k_0} \otimes v_1^{\otimes k_1}$}. Since $c_0\leq k_1$, then $x_{2\underline{1}}(0)^{c_0} v_\Lambda\neq 0$.
Moreover, $x_{2\underline{1}}(0) v_1 = v_2$ and $x_{2\underline{1}}(0) v_2 = 0$.

Hence,
\begin{equation} \label{trans}
x_{12}(0)^N x(\pi)v_\Lambda = \ox(\pi_1) x_{2\underline{1}}(0)^{c_0} v_\Lambda.
\end{equation}
Furthermore,
$x_{12}(0)^{N+1} x(\pi)v_\Lambda = 0$.

Different possibilities of distribution of $x_{2\underline{1}}(0)$'s on the tensor product \smash{$v_0^{\otimes k_0} \otimes v_1^{\otimes k_1}$} in \eqref{eq3} will be handled by certain coefficients of intertwining operators.

\subsection{Intertwining operators}

For \smash{$L_{C_2^{(1)}}(\Lambda_1)$} there is an operator (a coefficient of an intertwining operator) $w\colon \smash{L_{C_2^{(1)}}(\Lambda_1)}\to \smash{L_{C_2^{(1)}}(\Lambda_2)}$ such that
$v_1 \xrightarrow{w} 0$, $ v_2 \xrightarrow{w} v_{12}$,
where $v_{12}$ is a highest weight vector of \smash{$L_{C_2^{(1)}}(\Lambda_2)$}, and $w$ commutes with the action of $\gt_1$ (see \cite[Proposition~7]{BPT} or \cite[Remark~6.3]{P}).

On $L(\Lambda)$ use tensor products of these operators
\[w_{k_1,s}=\underbrace{1\otimes\cdots\otimes 1}_{k_0}\otimes \underbrace{1\otimes\cdots\otimes 1}_{k_1 - s}\otimes \underbrace{w\otimes\cdots\otimes w}_{s},\]
where $s\leq k_1$

Act on \eqref{eq3} by $w_{k_1,c_0}$
\begin{gather}
	w_{k_1,c_0}  \ox(\pi_1) x_{2\underline{1}}(0)^{c_0} v_\Lambda\nonumber\\
\qquad = \ox(\pi_1) w_{k_1,c_0} x_{2\underline{1}}(0)^{c_0} v_\Lambda\nonumber \\
	\qquad = \ox(\pi_1) w_{k_1,c_0} v_0^{\otimes k_0}\otimes v_1^{\otimes (k_1-c_0)}\otimes v_2^{\otimes c_0} + \sum \ox(\pi_1) w_{k_1,c_0} v_0^{\otimes k_0}\otimes \cdots \label{eq4}\\
	\qquad = \ox(\pi_1) \underbrace{v_0^{\otimes k_0}\otimes v_1^{\otimes (k_1-c_0)}\otimes v_{12}^{\otimes c_0}}_{v_{\Lambda'}},\nonumber
\end{gather}
where $\Lambda'=k_0 \Lambda_0+(k_1-c_0)\Lambda_1 + c_0 \Lambda_2$.
In \eqref{eq4}, the sum goes over all other distributions of~$x_{2\underline{1}}(0)$'s on tensor factors. These have at least one $v_1$ among the last $c_0$ tensor factors and hence are annihilated by $w_{k_1,c_0}$.

Since $x(\pi)$ satisfies initial conditions \eqref{ICA1}--\eqref{ICA2} for \smash{$L_{A_1^{(1)}}(\Lambda)$}, then,
by \eqref{DCA1.5} and \eqref{DCA2} for~${i=0}$,
$\ox(\pi_1)$ satisfies initial conditions \eqref{ICC1}--\eqref{ICC2} for \smash{$W_{C_2^{(1)}}(\Lambda')$}, with $k_0'=k_0$, $k_1'=k_1-c_0$, $k_2'=c_0$.

Hence
\smash{$\ox(\pi_1)x_{2\underline{1}}(0)^{c_0} v_\Lambda \xrightarrow{w_{k_1,c_0}} \ox(\pi_1) v_{\Lambda'}$}
and $\ox(\pi_1)$ satisfies initial conditions for \smash{$W_{C_2^{(1)}}(\Lambda')$}.
Note also that
\smash{$\ox(\pi_1)x_{2\underline{1}}(0)^{c_0} v_\Lambda \xrightarrow{w_{k_1,s}} 0$},
for $c_0 < s\leq k_1$.

\subsection{Proof of linear independence}
	
Let
\begin{equation} \label{eqLinIndep}
	\sum_\pi C_\pi x(\pi) v_\Lambda=0,
\end{equation}
be a relation of linear dependence, where all monomials in \eqref{eqLinIndep} satisfy difference and initial conditions for \smash{$L_{A_1^{(1)}}(\Lambda)$}.	
Let $ N=\max_\pi N(\pi)$.
Proceed inductively on $N$; act on \eqref{eqLinIndep} by~$x_{12}(0)^N$
\begin{align}
0 & = \sum_{N(\pi)<N} C_\pi x_{12}(0)^N x(\pi) v_\Lambda + \sum_{N(\pi)=N} C_\pi x_{12}(0)^N x(\pi) v_\Lambda
 = \sum_{N(\pi)=N} C_\pi x_{12}(0)^N x(\pi) v_\Lambda\nonumber\\
\label{eq5}	& = \sum_{\substack{N(\pi)=N\\ c_0=0}} C_\pi \ox(\pi_1) v_\Lambda +\sum_{\substack{N(\pi)=N\\ c_0=1}} C_\pi \ox(\pi) x_{2\underline{1}} (0) v_\Lambda 
+ \dots + \sum_{\substack{N(\pi)=N\\ c_0=s}} C_\pi \ox(\pi_1) x_{2\underline{1}} (0)^s v_\Lambda,
\end{align} 	
	for some $s\leq k_1$.
	The second equality follows since $x_{12}(0)^N x(\pi)=0$ if $N(\pi)<N$, while the third follows from \eqref{trans}.
	Proceed inductively on $s$; act on \eqref{eq5} by $w_{k_1,s}$. Then all sums except the last one are annihilated. One gets
\begin{equation} \label{eqLinIndep'}
\sum_{\substack{N(\pi)=N\\c_0=s}} C_\pi \ox(\pi_1) v_{\Lambda'}=0,
\end{equation}
 for $\Lambda'=k_0\Lambda_0+(k_1-s)\Lambda_1+ s\Lambda_2$.	
This is a relation of linear dependence in \smash{$W_{C_2^{(1)}}(\Lambda')$}. Note that all monomials in \eqref{eqLinIndep'} satisfy difference and initial conditions for \smash{$W_{C_2^{(1)}}(\Lambda')$}. Since these are linearly independent, all $C_\pi$ in \eqref{eqLinIndep'} are zero.

 \subsection*{Acknowledgments}
We are deeply honored to contribute this work to an issue dedicated to Jim Lepowsky -- an exceptional mathematician whose work has significantly shaped our scientific paths, a wonderful friend, and one of the kindest people that we know.

This work was supported by the project ``Implementation of cutting-edge research and its application as part of the Scientific Center of Excellence for Quantum and Complex Systems, and Representations of Lie Algebras'', PK.1.1.02, European Union, European Regional Development Fund. Also, this work has been supported by Croatian Science Foundation under the project IP-2022-10-9006.

\pdfbookmark[1]{References}{ref}
\LastPageEnding


\begin{thebibliography}{99}
\footnotesize\itemsep=0pt

\bibitem{BPT}
Baranovi\'c I., Primc M., Trup\v{c}evi\'c G., Bases of
  {F}eigin--{S}toyanovsky's type subspaces for~{$C_\ell^{(1)}$},
  \href{https://doi.org/10.1007/s11139-016-9840-y}{\textit{Ramanujan~J.}}
  \textbf{45} (2018), 265--289,
  \href{http://arxiv.org/abs/1603.04594}{arXiv:1603.04594}.

\bibitem{CMPP}
Capparelli S., Meurman A., Primc A., Primc M., New partition identities
  from~{$C^{(1)}_\ell$}-modules,
  \href{https://doi.org/10.3336/gm.57.2.01}{\textit{Glas. Mat. Ser.}}
  \textbf{57} (2022), 161--184,
  \href{http://arxiv.org/abs/2106.06262}{arXiv:2106.06262}.

\bibitem{FKLMM}
Feigin B., Kedem R., Loktev S., Miwa T., Mukhin E., Combinatorics of
  the~{$\widehat{\mathfrak{sl}}_2$} spaces of coinvariants,
  \href{https://doi.org/10.1007/BF01236061}{\textit{Transform. Groups}}
  \textbf{6} (2001), 25--52,
  \href{http://arxiv.org/abs/math-ph/9908003}{arXiv:math-ph/9908003}.

\bibitem{K}
Kac V.G., Infinite-dimensional {L}ie algebras, 3rd~ed.,
  \href{https://doi.org/10.1017/CBO9780511626234}{Cambridge University Press}, Cambridge, 1990.

\bibitem{MP}
Meurman A., Primc M., Annihilating fields of standard modules
  of~{$\mathfrak{sl}(2,\mathbb{C})^\sim$} and combinatorial identities,
  \href{https://doi.org/10.1090/memo/0652}{\textit{Mem. Amer. Math. Soc.}}
  \textbf{137} (1999), viii+89~pages,
  \href{http://arxiv.org/abs/math.QA/9806105}{arXiv:math.QA/9806105}.

\bibitem{P}
Primc M., Combinatorial bases of modules for affine {L}ie
  algebra~{$B_2^{(1)}$},
  \href{https://doi.org/10.2478/s11533-012-0111-x}{\textit{Cent. Eur.~J.
  Math.}} \textbf{11} (2013), 197--225,
  \href{http://arxiv.org/abs/1002.3535}{arXiv:1002.3535}.

\bibitem{PT}
Primc M., Trup\v{c}evi\'c G., Linear independence for~{$C_\ell^{(1)}$} by
  using~{$C_{2\ell}^{(1)}$},
  \href{https://doi.org/10.1016/j.jalgebra.2024.08.003}{\textit{J.~Algebra}}
  \textbf{661} (2025), 341--356,
  \href{http://arxiv.org/abs/2403.06881}{arXiv:2403.06881}.

\end{thebibliography}
\end{document}